\documentclass{jnmp}
\usepackage{amsmath}

\JNMPnumberwithin{equation}{section}

\newtheorem{theorem}{Theorem}
\newtheorem{lemma}{Lemma}
\theoremstyle{definition}
\newtheorem*{remark}{Remark}

\begin{document}

\renewcommand{\evenhead}{M Liefvendahl and G Kreiss}

\renewcommand{\oddhead}{Bounds for the Threshold Amplitude for Plane
Couette Flow}

\thispagestyle{empty}

\FirstPageHead{9}{3}{2002}{\pageref{Liefvendahl-firstpage}--\pageref{Liefvendahl-lastpage}}{Article}

\copyrightnote{2002}{M Liefvendahl and G Kreiss}

\Name{Bounds for the Threshold Amplitude\\ for Plane Couette Flow}
\label{Liefvendahl-firstpage}

\Author{Mattias LIEFVENDAHL and Gunilla KREISS}

\Address{Department of Numerical Analysis and Computing
Science, \\ KTH, Stockholm S-100 44, Sweden \\
E-mail: mli@nada.kth.se, \  gunillak@nada.kth.se}

\Date{Received October 3, 2001; Accepted April 15, 2002}

\begin{abstract}
\noindent
We prove nonlinear stability for finite amplitude perturbations
of plane Couette flow. A bound of the solution of the resolvent
equation in the unstable complex half-plane is used to estimate
the solution of the full nonlinear problem.
The result is a lower bound, including Reynolds number dependence,
of the threshold amplitude below
which all perturbations are stable. Our result is an
improvement of the corresponding bound derived in~\cite{KLH:1}.
\end{abstract}

\section{Introduction}

Plane Couette flow is stable to infinitesimal perturbations
for all Reynolds numbers~\cite{Romanov:1}.
Finite amplitude perturbations on the other hand can induce
transition to turbulence. Thus there is a threshold amplitude
below which all perturbations decay eventually (there may be
transient growth) and above which some perturbations lead
to turbulent flow. In the present paper we derive a bound
depending on the Reynolds number~$R$, for this threshold.

The stability of shear flows, and plane Couette flow in particular,
has been extensively studied, we refer to \cite{Lin:1,DR:1}
and \cite{SH:1} for more background information.
Early work focused on the eigenvalues of the linearized problem.
The question of the asymptotic value of the threshold amplitude
for large~$R$ was starting to be investigated in the early
1990s~\cite{TTRD:1}. A~dependence $\mathcal{O}(R^{-\gamma})$
with $\gamma\approx 5/4$ is supported by computations in~\cite{RSBH:1}.
The asymptotic analysis of~\cite{Chapman:1} gives the
result~$\gamma\approx 1$. Another approach to
this question is found in~\cite{KLH:1} where an estimate of
the resolvent equation is used to prove an upper bound for
the exponent, $\gamma\leq 21/4$. In this paper we use the same approach
and sharpen the bound of~$\gamma$ obtained in~$\cite{KLH:1}$.
This is possible because we work with a norm which weighs
the different coefficients of the velocity vector, the weights
depend on the Reynolds number. It is numerically demonstrated
in~\cite{jag:2} that, in the weighted norm, the resolvent grows
linearly with~$R$, as opposed to quadratically which is
the case of the energy norm used in~\cite{KLH:1}.
Because of the weighting the exponent $\gamma$ will be different
for the different components of the velocity vector. The precise
result is stated in~Theorem \ref{thm:1}.

The proof of the result of this paper is carried out
with the same techniques as
in~\cite{KLH:1}. However, the starting point is a different resolvent
estimate. Also, the energy norm fits the structure of the Navier--Stokes
equations, the modified norm we use does not. This leads to several
technical complications.

In Section \ref{sec:2} we give the mathematical formulation of
the problem, the Navier-Stokes equations, the linearized equation
and the resolvent equation. Then we state the main result in Theorem~\ref{thm:1}.
The proof of the theorem is organized as follows.
In Section~\ref{sec:3} we analyze the linear problem.
Using the result of~\cite{jag:2} mentioned above, we
derive an estimate for higher derivatives of
the solution of the resolvent equation
(this solution is the Laplace transform of the velocity field).
By inverse transformation we obtain a time dependent estimate.
In Section~\ref{sec:4} the linear estimate is used to derive the
threshold value for the nonlinear problem.

In Appendix~\ref{app:1} we have collected the definition
of the weighted norm and other non-standard norms used
in the paper.
Technical results
including estimates of norms of the nonlinearity are collected
in Appendix~\ref{app:2}.
In Appendix~\ref{app:3}, we prove a theorem which will
be applied to show that the perturbation tends to zero if the
suitable a priori estimate holds. The result of Appendix~\ref{app:3}
explains the choice of left hand side in the linear estimate
of Theorem~\ref{thm:2} in Section~\ref{sec:3}.

\section{The main result}\label{sec:2}

We choose the coordinate system so that the velocity field of
Couette flow is given by
\[
\boldsymbol{u}_{\mbox{\rm \tiny Co}}=\begin{pmatrix} x_2\\ 0\\
0
\end{pmatrix}
\]
in the domain
\[
\Omega=\left\{\boldsymbol{x}\in\mathbb{R}^3: \ -1<x_2<1\right\}.
\]
We use bold letters to denote vectors and subscripts to identify
the components. The functions $\boldsymbol{u}_{\mbox{\rm \tiny
Co}}$ and $p_{\mbox{\rm \tiny Co}}={\rm const}$ constitute a
stationary solution of the nondimensionalized Navier--Stokes
equations
\begin{gather}
\boldsymbol{u}_t+\boldsymbol{G}(\boldsymbol{u})+
{\rm grad}\, p  =  \frac{1}{R}\Delta \boldsymbol{u},\nonumber\\
 \label{eq:m5}
{\rm div}\, \boldsymbol{u} =  0.
\end{gather}
Here $R$ denotes the Reynolds number and $\boldsymbol{G}$ is the following
nonlinear differential operator.
\[
\boldsymbol{G}(\boldsymbol{u})=\sum_{k=1}^3 u_k \frac{\partial \boldsymbol{u}}{\partial x_k}.
\]
We have no-slip boundary conditions, $\boldsymbol{u}=0$ on $\partial \Omega$,
and the initial condition
\[
\boldsymbol{u}(\boldsymbol{x},0)=\boldsymbol{u}_{\mbox{\rm \tiny
Co}}(\boldsymbol{x})+ \boldsymbol{v}^{(0)}(\boldsymbol{x}),
\]
where $\boldsymbol{v}^{(0)}$ is the initial perturbation.
We will assume that $\boldsymbol{v}^{(0)}\in H^4$. Local existence
of a classical solution of~(\ref{eq:m5}) is proven in~\cite{Ladyzhenskaya:1} (theorem
5 together with theorem 7 on p.~161 and p.~167 respectively of~\cite{Ladyzhenskaya:1}).
Furthermore we will derive a priori estimates which (for sufficiently small
initial data) allows extension of the local solution to a global solution
by successive application of the above result.

We investigate the stability of Couette flow by linearization of
the Navier--Stokes equations. We denote the perturbation by $\boldsymbol{v}$
which thus is related to $\boldsymbol{u}$ according to
\[
\boldsymbol{u}=\boldsymbol{u}_{\mbox{\rm \tiny
Co}}+\boldsymbol{v}.
\]
Below we will apply the Laplace transform and for this purpose we
want a problem with homogeneous initial data. This is accomplished
with the introduction of $\boldsymbol{w}$ by
\begin{equation} \label{eq:m10}
\boldsymbol{v}=e^{-t}\boldsymbol{v}^{(0)}+\boldsymbol{w}.
\end{equation}
The function $\boldsymbol{w}$ satisfies the following problem
\begin{gather}
\boldsymbol{w}_t+{\rm grad}\, p  =  L\boldsymbol{w} +e^{-t}(L+1)\boldsymbol{v}^{(0)}-
\boldsymbol{G}\left(e^{-t}\boldsymbol{v}^{(0)}+\boldsymbol{w}\right),  \nonumber\\
{\rm div}\, \boldsymbol{w}  =  0, \nonumber \\
\boldsymbol{w}(\boldsymbol{x},0)  =  0,  \label{eq:m2}
\end{gather}
where
\[
L=\frac{1}{R}\Delta -x_2\frac{\partial }{\partial x_1}-\left(\begin{array}{ccc}
0 & 1 & 0 \\ 0 & 0 & 0 \\ 0 & 0 & 0
\end{array}\right).
\]
For $\boldsymbol{w}$ we also have no-slip boundary conditions.
If we assume that $\boldsymbol{v}^{(0)}$ and $\boldsymbol{w}$ are small and neglect
$\boldsymbol{G}$ in~(\ref{eq:m2}) then we see
that~$\boldsymbol{w}$ satisfies a linear problem with forcing depending
on the initial perturbation. This is the linearized equation
for Couette flow.

We will first consider the forcing in (\ref{eq:m2}) as a given function,
denoted by $\boldsymbol{f}=\boldsymbol{f}(\boldsymbol{x},t)$ and derive estimates for the
solution. Later we will return to the particular form of the forcing
in~(\ref{eq:m2}) and its dependence on the initial data and~$\boldsymbol{w}$.

Application of the Laplace transform now give the resolvent equation
\begin{gather}
s\hat{\boldsymbol{w}}+{\rm grad}\, \hat{p}  =
 L\hat{\boldsymbol{w}}+\hat{\boldsymbol{f}},\nonumber \\
{\rm div}\, \hat{\boldsymbol{w}}  =   0,  \nonumber\\
\hat{\boldsymbol{w}}  = 0,\qquad \boldsymbol{x}\in \partial \Omega. \label{eq:m3}
\end{gather}
This problem was investigated in \cite{jag:2} where
the following estimate of the solution was obtained
\begin{equation} \label{eq:m1}
\|\hat{\boldsymbol{w}}\|_{m}^2 \leq C R^2 \|\hat{\boldsymbol{f}}\|_m^2,
\qquad  {\rm Re}\, s\geq 0.
\end{equation}
See Appendix \ref{app:1} equation (\ref{eq:a4}) for the definition of
the modified ($m$-)norm.
We have collected the definitions of all norms used in this paper
in Appendix~\ref{app:1}.
Throughout the paper we will use $C$ to denote
constants which appear in inequalities, we do not use subscripts
to identify different constants. We emphasize that $C$ denotes an ``absolute
constant'' and does not depend on the forcing, the initial data
or the Reynolds number, which is the only
parameter in the problem.

We now state our main result in the following theorem.

\begin{theorem} \label{thm:1}
The resolvent estimate (\ref{eq:m1}) implies that there is
a constant $\delta>0$ such that
\[
\lim_{t\rightarrow \infty}\|\boldsymbol{v}(\cdot,t)\|_{\infty}=0
\]
if
\begin{gather*}
\|v_k(\cdot,0)\|_{H^4}  \leq  \frac{\delta}{R^3},\qquad k=1,3, \\
\|v_2(\cdot,0)\|_{H^4}  \leq  \frac{\delta}{R^4}.
\end{gather*}
\end{theorem}

Comparing this theorem to the previously obtained bound on the
threshold in~\cite{KLH:1} we see that the $R$-exponent for the
first and third $\boldsymbol{v}$-components is improved from 5.25
(in~\cite{KLH:1}) to~3 and the exponent for the second component
is improved from 5.25 to~4.

It may appear excessive to require that the initial data is as
smooth as $H^4$. The result in~\cite{Romanov:1} require the
$H^2$-norm of the initial data to be small. However,
in~\cite{Romanov:1} a completely different approach is used with
the application of semi-group methods and we are not aware of any
results on the threshold amplitude obtained with this approach. It
may be possible to weaken our smoothness requirement using an
appropriate local existence theorem for the Navier--Stokes
equations which incorporates the smoothing property. We also note
the error in \cite[p.~193]{KLH:1} where it is stated that small
$H^2$-norm implies stability, for their result to hold small
$H^4$-norm is required.

If one is interested in
a particular perturbation $\boldsymbol{v}^{(0)}(\boldsymbol{x})
=\alpha \boldsymbol{\varphi}(\boldsymbol{x})$,
and the threshold coefficient $\alpha(R)$, then it is insignificant
which $H^p$-norm is required since
\[
\big\|\boldsymbol{v}^{(0)}\big\|_{H^p} \leq \frac{1}{R^\gamma}
\qquad \Rightarrow \qquad \alpha(R)\leq \frac{1}{R^\gamma
\|\boldsymbol{\varphi}\|_{H^p}}= \mathcal{O}(R^{-\gamma}).
\]

\section{The linear estimate} \label{sec:3}

In this section we derive estimates for the solution of the
following problem in terms of the forcing function
$\boldsymbol{f}$
\begin{gather}
\boldsymbol{w}_t+{\rm grad}\, p  =  L\boldsymbol{w} +\boldsymbol{f},\nonumber  \\
{\rm div}\, \boldsymbol{w}  =  0, \nonumber \\
\boldsymbol{w}(\boldsymbol{x},0)  =  0. \label{eq:l1}
\end{gather}
The results are collected in the inequality of the theorem below.
In the statement of the theorem the $\tilde{H}_1$- and
$\tilde{H}_2$-norms appear for the first time, their definition is
given in Appendix~\ref{app:1}.

\begin{theorem} \label{thm:2}
The resolvent estimate~(\ref{eq:m1}) implies that the solution of
equation~(\ref{eq:l1}) satisfies
\begin{gather}
\|\boldsymbol{w}(\cdot,T)\|_{\tilde{H}_2}^2+\int_0^T\left(
\|\boldsymbol{w}(\cdot,t)\|_{\tilde{H}_1}^2+
\|\boldsymbol{w}_t(\cdot,t)\|_{\tilde{H}_1}^2\right)dt \nonumber\\
\phantom{\|\boldsymbol{w}(\cdot,T)\|_{\tilde{H}_2}^2}{}\leq
C\left[
R\big\|\boldsymbol{f}^{(0)}\big\|_{\tilde{H}_2}^2+R^2\|\boldsymbol{f}(\cdot,T)\|^2+
R^2\big\|(L+1)\boldsymbol{f}^{(0)}\big\|_m^2\right] \nonumber \\
\phantom{\|\boldsymbol{w}(\cdot,T)\|_{\tilde{H}_2}^2}{}+
CR^2\int_0^T\left(
\|\boldsymbol{f}(\cdot,t)\|_{m}^2+\|\boldsymbol{f}_t(\cdot,t)\|_{m}^2\right)dt.
\label{eq:l21}
\end{gather}
\end{theorem}

For the remainder of the paper we will suppress the notation
``$(\cdot,t)$'' which indicates at which
time a norm is evaluated, this will be clear from the context.

In Section~\ref{sec:31} we analyze the Laplace transformed
problem~(\ref{eq:m3}) and then in Section~\ref{sec:32} we use
these results to conclude the proof of Theorem~\ref{thm:2}.

\subsection{Estimates of the transformed functions} \label{sec:31}

The starting point for this section is inequality (\ref{eq:m1}) which
we will use to derive the following lemma.

\begin{lemma}
If ${\rm Re}\, s\geq 0$, the solution of (\ref{eq:m3}) satisfies
\begin{equation} \label{eq:l15}
\|\hat{\boldsymbol{w}}\|_{\tilde{H}_1}^2\leq
CR^2\|\hat{\boldsymbol{f}}\|_m^2.
\end{equation}
The estimate is uniform in $s$.
\end{lemma}

\begin{proof} By partial integration, the following identity
is easily derived
\begin{equation} \label{eq:l2}
{\rm Re}\,
(\hat{\boldsymbol{w}},L\hat{\boldsymbol{w}})=-\frac{1}{R}\sum_{k=1}^3
\left\|\frac{\partial \hat{\boldsymbol{w}}}{\partial
x_k}\right\|^2+{\rm Re}\,(\hat{w}_1,\hat{w}_2).
\end{equation}
Now, we take the $L^2$-inner product of $\hat{\boldsymbol{w}}$ and
the resolvent equation~(\ref{eq:m3}). Using~(\ref{eq:l2}) and the
fact that~$\hat{\boldsymbol{w}}$ is solenoidal we obtain
\begin{equation} \label{eq:l3}
{\rm Re}\, s\,\|\hat{\boldsymbol{w}}\|^2+\frac{1}{R}\sum_{k=1}^3
\left\|\frac{\partial \hat{\boldsymbol{w}}}{\partial
x_k}\right\|^2= {\rm Re} \left[
(\hat{w}_1,\hat{w}_2)+(\hat{\boldsymbol{w}},\hat{\boldsymbol{f}})
\right].
\end{equation}
The right hand side in~(\ref{eq:l3}) can be estimated using the
inequality $ab\leq \epsilon a^2/2+b^2/2\epsilon$ which holds for
all $\epsilon>0$. Without explicit mention we will use this
inequality several times below. Here we take
$a=\|\hat{\boldsymbol{w}}\|$, $b=\|\hat{\boldsymbol{f}}\|$ and
$\epsilon=1/R$. We also use the fact that the $L_2$-norm is
smaller than the modified norm ($R\geq 1$), and obtain
\begin{equation} \label{eq:l4}
{\rm Re} \left[(\hat{w}_1,\hat{w}_2)
+(\hat{\boldsymbol{w}},\hat{\boldsymbol{f}}) \right] \leq
CR\|\hat{\boldsymbol{f}}\|_m^2.
\end{equation}
Combining (\ref{eq:l4}) with (\ref{eq:l3}) and using ${\rm Re}\,
s\geq 0$ we have derived
\begin{equation} \label{eq:l12}
\sum_{k=1}^3 \left\|\frac{\partial \hat{\boldsymbol{w}}}{\partial
x_k}\right\|^2\leq CR^2\|\hat{\boldsymbol{f}}\|_m^2.
\end{equation}
It remains to estimate $J_2(\hat{\boldsymbol{w}})$ (defined in
Appendix~\ref{app:2}). To do this we differentiate the resolvent
equation with respect to $x_1$ and $x_3$. We give the details in
the $x_1$-case. We obtain
\begin{equation} \label{eq:l9}
s\frac{\partial \hat{\boldsymbol{w}}}{\partial x_1}+{\rm grad}
\frac{\partial \hat{p}}{\partial x_1}=L\frac{\partial
\hat{\boldsymbol{w}}}{\partial x_1}+ \frac{\partial
\hat{\boldsymbol{f}}}{\partial x_1}.
\end{equation}
Note that $L$ commutes with $\partial/\partial x_1$ and
$\partial/\partial x_3$ but not $\partial/\partial x_2$. We take
the inner product of~(\ref{eq:l9}) with
$\partial\hat{\boldsymbol{w}}/\partial x_1$ and obtain
\begin{equation} \label{eq:l11}
{\rm Re}\, s\left\|\frac{\partial \hat{\boldsymbol{w}}}{\partial
x_1}\right\|^2+
\frac{1}{R}\sum_{k=1}^3\left\|\frac{\partial^2\hat{\boldsymbol{w}}}{\partial
x_1 \partial x_k}\right\|^2= {\rm Re}\left[
\left(\frac{\partial \hat{w}_1}{\partial x_1},\frac{\partial
\hat{w}_2}{\partial x_1}\right)+ \left(\frac{\partial
\hat{\boldsymbol{w}}}{\partial x_1},\frac{\partial
\hat{\boldsymbol{f}}}{\partial x_1}\right) \right].
\end{equation}
The right hand side is estimated in the following way.
\begin{gather}
{\rm Re}\left[ \left(\frac{\partial \hat{w}_1}{\partial
x_1},\frac{\partial \hat{w}_2}{\partial x_1}\right)+
\left(\frac{\partial \hat{\boldsymbol{w}}}{\partial
x_1},\frac{\partial \hat{\boldsymbol{f}}}{\partial x_1}\right)
\right] \nonumber\\
\qquad {}\leq
\frac{1}{4R}\left\|\frac{\partial^2\hat{w}_1}{\partial
x_1^2}\right\|^2+ R\|\hat{w}_2\|^2+
\frac{1}{4R}\left\|\frac{\partial^2\hat{\boldsymbol{w}}}{\partial
x_1^2}\right\|^2+ R\|\hat{\boldsymbol{f}}\|^2. \label{eq:l10}
\end{gather}
Now we insert the estimate (\ref{eq:l10}) into (\ref{eq:l11})
and cancel the second
derivatives in the right hand side with the corresponding terms
on the left hand side. Further simplification leads to
\[
\frac{1}{R}\sum_{k=1}^3\left\|\frac{\partial^2\hat{\boldsymbol{w}}}{\partial
x_1 \partial x_k}\right\|^2\leq
C\left(\frac{1}{R}\|\hat{\boldsymbol{w}}\|_{m}^2+R\|\hat{\boldsymbol{f}}\|_m^2\right)\leq
CR\|\hat{\boldsymbol{f}}\|_m^2.
\]
The same procedure applies to the resolvent equation differentiated
with respect to $x_3$. We have proven
\begin{equation} \label{eq:l13}
J_2(\hat{\boldsymbol{w}}) \leq CR^2\|\hat{\boldsymbol{f}}\|_m^2.
\end{equation}
A combination of (\ref{eq:m1}), (\ref{eq:l12}) and (\ref{eq:l13})
gives the lemma.
\end{proof}

\begin{remark} To obtain (\ref{eq:l2}) we used
\begin{equation} \label{eq:l30}
(\hat{\boldsymbol{w}},{\rm grad}\, \hat{p})=-({\rm div}\,
\hat{\boldsymbol{w}},\hat{p})=0,
\end{equation}
where the first equality follows by partial integration and the
second since $\hat{\boldsymbol{w}}$ is solenoidal (${\rm
div}\,\hat{\boldsymbol{w}}=0$). It would be preferable to work
with modified ($m$-)norms throughout the paper, the reason one
cannot do this is~(\ref{eq:l30}). If one takes the modified inner
product of~$\hat{\boldsymbol{w}}$ and the resolvent equation then
it is not possible to eliminate the pressure using~(\ref{eq:l30}).
The modified inner product is defined by
\[
(\boldsymbol{u},\boldsymbol{v})_m=(u_1,v_1)+R^2(u_2,v_2)+(u_3,v_3).
\]
This is the explanation of the simplifications with the use of the energy
($L_2$) norm mentioned in the introduction.
\end{remark}

\subsection{Estimates of the time dependent functions} \label{sec:32}

Before we prove Lemma~\ref{lemma:2} below and Theorem~\ref{thm:2}
we state three inequalities which will be of use. The following is
derived from equation~(\ref{eq:l1}) using~(\ref{eq:l2})
\begin{equation} \label{eq:l5}
\frac{d}{dt}\|\boldsymbol{w}\|^2+\frac{1}{R}\sum_{k=1}^3
\left\|\frac{\partial \boldsymbol{w}}{\partial x_k}\right\|^2 \leq
|(\boldsymbol{w},\boldsymbol{f})|+|(w_1,w_2)|.
\end{equation}
In the same way we obtain the next inequality, for $l=1,3$,
starting from~(\ref{eq:l1}) differentiated with respect to $x_1$
and~$x_3$ respectively
\begin{equation} \label{eq:l6}
\frac{d}{dt} \left\|\frac{\partial \boldsymbol{w}}{\partial
x_l}\right\|^2+\frac{1}{R}\sum_{k=1}^3
\left\|\frac{\partial^2\boldsymbol{w}}{\partial x_l
\partial x_k}\right\|^2 \leq
\left|\left(\frac{\partial^2\boldsymbol{w}}{\partial
x_l^2},\boldsymbol{f} \right)\right|+ \left|\left(\frac{\partial^2
w_1}{\partial x_l^2},w_2 \right)\right|.
\end{equation}
We apply Plancherel's formula for the Laplace
transform to (\ref{eq:l15}), with
the imaginary axis as integration contour.
This yields
\begin{equation} \label{eq:l7}
\int_0^T \|\boldsymbol{w}\|_{\tilde{H}_1}^2 dt \leq CR^2\int_0^T
\|\boldsymbol{f}\|_{m}^2dt,
\end{equation}
which holds for all $T>0$, see \cite[p.~235--239]{KL:1}  for this
application of Plancherel's formula.

\begin{lemma} \label{lemma:2}
The solution of (\ref{eq:l1}) satisfies
\begin{equation} \label{eq:l16}
\left\|\boldsymbol{w}\right\|^2+\left\|\frac{\partial
\boldsymbol{w}}{\partial x_1}\right\|^2+ \left\|\frac{\partial
\boldsymbol{w}}{\partial x_3}\right\|^2 \leq CR \int_0^T
\|\boldsymbol{f}\|_{m}^2dt,
\end{equation}
for all $T>0$.
\end{lemma}

\begin{proof}
From (\ref{eq:l5}) we obtain
\[
\frac{d}{dt}\|\boldsymbol{w}\|^2\leq
\frac{\epsilon}{2}\|\boldsymbol{w}\|_{m}^2+
\frac{1}{2\epsilon}\|\boldsymbol{f}\|_{m}^2+\frac{1}{R}\|\boldsymbol{w}\|_{m}^2,
\qquad \epsilon>0.
\]
We integrate this from $t=0$ to $T$, choose $\epsilon=1/R$ and use (\ref{eq:l7}),
\[
\|\boldsymbol{w}\|^2 \leq CR\int_0^T \|\boldsymbol{f}\|_{m}^2dt.
\]
We estimate the right hand side in (\ref{eq:l6}) according to
\[
\left|\left(\frac{\partial^2 w_1}{\partial x_l^2},w_2
\right)\right|+ \left|\left(\frac{\partial^2
\boldsymbol{w}}{\partial x_l^2},\boldsymbol{f} \right)\right|\leq
\frac{1}{4R}\left\|\frac{\partial^2 w_1}{\partial
x_l^2}\right\|^2+ R\|w_2\|^2+\frac{1}{4R}\left\|\frac{\partial^2
\boldsymbol{w}}{\partial x_l^2}\right\|^2+ R\|\boldsymbol{f}\|^2,
\]
recall that $l=1,3$. This inequality is inserted
into~(\ref{eq:l6}). We move the second derivative terms to the
left hand side and integrate in time to obtain
\[
\left\|\frac{\partial \boldsymbol{w}}{\partial x_l}\right\|^2\leq
CR\int_0^T \|\boldsymbol{f}\|_{m}^2dt.
\]
This completes the proof of Lemma~\ref{lemma:2}.
\end{proof}

Now we turn to the proof of Theorem~\ref{thm:2}. The function
$\boldsymbol{w}_t$ satisfies the same PDE as $\boldsymbol{w}$
(with $p$ and $\boldsymbol{f}$ replaced with $p_t$ and
$\boldsymbol{f}_t$ respectively) and inhomogeneous initial data,
\[
\boldsymbol{w}_t(\boldsymbol{x},0)=
 \boldsymbol{f} (\boldsymbol{x},0)=:\boldsymbol{f}^{(0)}(\boldsymbol{x}).
\]
We use the symbol $:=$ to indicate that the expression on the ``colon side''
is defined in terms of the expression on the other side. Now we introduce
$\boldsymbol{\omega}$ by
\begin{equation} \label{eq:l17}
\boldsymbol{w}_t=e^{-t}\boldsymbol{f}^{(0)}+\boldsymbol{\omega}.
\end{equation}
The function $\boldsymbol{\omega}$ satisfies the following problem.
\begin{gather*}
\boldsymbol{\omega}_t+{\rm grad}\, p_t  =
 L\boldsymbol{\omega} +\boldsymbol{f}_t+e^{-t}(L+1)\boldsymbol{f}^{(0)}, \\
{\rm div}\, \boldsymbol{\omega}  =  0, \\
\boldsymbol{\omega}(\boldsymbol{x},0)  =  0.
\end{gather*}
The only difference between this problem and equation
(\ref{eq:l1}) for $\boldsymbol{w}$ is the forcing. The
inequalities~(\ref{eq:l7}) and~(\ref{eq:l16}) can thus be applied
here yielding estimates for~$\boldsymbol{\omega}$.
Using~(\ref{eq:l17}) and the triangle inequality these can be
converted to the following estimates for~$\boldsymbol{w}$.
\begin{equation} \label{eq:l18}
\int_0^T \|\boldsymbol{w}_t\|_{\tilde{H}_1}^2dt \leq C\left(
\big\|\boldsymbol{f}^{(0)}\big\|_{\tilde{H}_1}^2+R^2\big\|(L+1)\boldsymbol{f}^{(0)}\big\|_m^2+
R^2\int_0^T \|\boldsymbol{f}_t\|_{m}^2 dt \right)
\end{equation}
and
\begin{gather}
\|\boldsymbol{w}_t\|^2+\left\|\frac{\partial
\boldsymbol{w}_t}{\partial x_l}\right\|^2 \nonumber\\
{}\leq  C \left(
\big\|\boldsymbol{f}^{(0)}\big\|^2+\left\|\frac{\partial
\boldsymbol{f}^{(0)}}{\partial x_l}\right\|^2
+R\big\|(L+1)\boldsymbol{f}^{(0)}\big\|_m^2+R\int_0^T
\|\boldsymbol{f}_t\|_{m}^2dt \right),\quad \!
l=1,3.\!\!\!\label{eq:l19}
\end{gather}
With (\ref{eq:l7}) and (\ref{eq:l18}) we have estimated the time
integral in the left hand side of the inequality of
Theorem~\ref{thm:2}. The $\tilde{H}_2$-norm at time $T$ remains.
To estimate the first order derivatives in the $\tilde{H}_2$-norm
we start from the following estimate which is easily obtained
from~(\ref{eq:l5})
\begin{equation} \label{eq:l20}
\sum_{k=1}^3\left\|\frac{\partial \boldsymbol{w}}{\partial
x_k}\right\|^2 \leq CR\left(
|(\boldsymbol{w},\boldsymbol{f})|+|(\boldsymbol{w},\boldsymbol{w}_t)|+|(w_1,w_2)|
\right).
\end{equation}
Using the Cauchy--Schwartz inequality and (\ref{eq:l16}) and
(\ref{eq:l19}) we see that the right hand side of~(\ref{eq:l20})
can be estimated in terms of the right hand side
of~(\ref{eq:l21}).

To estimate $J_2(\boldsymbol{w})$ we rearrange (\ref{eq:l6}) and obtain
\begin{equation} \label{eq:l22}
\sum_{k=1}^3 \left\|\frac{\partial^2\boldsymbol{w}}{\partial x_l
\partial x_k}\right\|^2 \leq CR \left[ \left|\left(
\frac{\partial w_1}{\partial x_l},\frac{\partial w_2}{\partial
x_l} \right) \right|+ \left|\left(
\frac{\partial^2\boldsymbol{w}}{\partial x_l^2},\boldsymbol{f}
\right) \right|+ \left|\left( \frac{\partial
\boldsymbol{w}}{\partial x_l},\frac{\partial
\boldsymbol{w}_t}{\partial x_l} \right) \right| \right],
\end{equation}
for $l=1,3$. For the second term in the right hand side
we have
\begin{equation} \label{eq:l23}
\left|\left( \frac{\partial^2 \boldsymbol{w}}{\partial
x_l^2},\boldsymbol{f} \right) \right| \leq
\frac{1}{4R}\left\|\frac{\partial^2 \boldsymbol{w}}{\partial
x_l^2}\right\|^2+R\|\boldsymbol{f}\|^2.
\end{equation}
We insert (\ref{eq:l23}) into (\ref{eq:l22}) and move the second
derivative term to the left hand side. The remaining two terms in
the right hand side of (\ref{eq:l22}) are estimated using the
Cauchy--Schwartz inequality and (\ref{eq:l16}) and (\ref{eq:l19}).
We see that $J_2(\boldsymbol{w})$ can be estimated by the
expression on the right hand side of the inequality of
Theorem~\ref{thm:2}. This concludes the proof of the theorem.

\section{Derivation of the threshold bound} \label{sec:4}

Now we return to the full nonlinear problem (\ref{eq:m2}) for $\boldsymbol{w}$.
This is the same as (\ref{eq:l1}) if we take
\begin{gather*}
\boldsymbol{f}=e^{-t}(L+1)\boldsymbol{v}^{(0)} \\
\phantom{\boldsymbol{f}=}{}+ \sum_{k=1}^3 \left(
e^{-2t}v^{(0)}_k\frac{\partial \boldsymbol{v}^{(0)}}{\partial
x_k}+ e^{-t}v^{(0)}_k\frac{\partial \boldsymbol{w}}{\partial x_k}+
e^{-t}w_k\frac{\partial \boldsymbol{v}^{(0)}}{\partial x_k}+
w_k\frac{\partial \boldsymbol{w}}{\partial x_k} \right).
\end{gather*}
For simplicity we first consider
\begin{equation} \label{eq:d1}
\boldsymbol{f}=e^{-t}(L+1)\boldsymbol{v}^{(0)}+\boldsymbol{G}(\boldsymbol{w}).
\end{equation}
The quadratic terms in $\boldsymbol{v}^{(0)}$ and the coupling
terms between $\boldsymbol{v}^{(0)}$ and $\boldsymbol{w}$ only add
technical difficulties, the asymptotic value of the threshold
amplitude is the same. To make the formulas below more lucid we
introduce the following notation for expressions in the inequality
of Theorem~\ref{thm:2}
\begin{gather*}
\mathcal{N}(\boldsymbol{w},T)  =
\|\boldsymbol{w}\|_{\tilde{H}_2}^2+\int_0^T
\left(\|\boldsymbol{w}\|_{\tilde{H}_1}^2+\|\boldsymbol{w}_t\|_{\tilde{H}_1}^2\right)dt, \\
\mathcal{M}(\boldsymbol{f},T)  = CR^2\left[\|\boldsymbol{f}\|^2
+\int_0^T\left(
\|\boldsymbol{f}\|_{m}^2+\|\boldsymbol{f}_t\|_{m}^2\right)dt\right],\\
\mathcal{M}_0(\boldsymbol{f},T)  = \mathcal{M}(\boldsymbol{f},T)+
C\left( R
\big\|\boldsymbol{f}^{(0)}\big\|_{\tilde{H}_2}^2+R^2\big\|(L+1)\boldsymbol{f}^{(0)}\big\|_m^2
\right).
\end{gather*}
As long as the solution $\boldsymbol{w}$ exists, the estimate of
Theorem~\ref{thm:2} is valid with forcing~(\ref{eq:d1}). Expressed
with our new notation we have
\begin{equation} \label{eq:d2}
\mathcal{N}(\boldsymbol{w},T) \leq
\mathcal{M}_0\left(e^{-t}(L+1)\boldsymbol{v}^{(0)},T\right)
+\mathcal{M}(\boldsymbol{G}(\boldsymbol{w}),T).
\end{equation}
Here we have used that $\boldsymbol{w}(\boldsymbol{x},0)=0$. The
first term on the right hand side of (\ref{eq:d2}) can be
estimated in terms of $\boldsymbol{v}^{(0)}$, independently of
$T$, as we show in Appendix~\ref{app:2},
inequality~(\ref{eq:a10}). We have
\begin{equation} \label{eq:d12}
\mathcal{M}_0\left(e^{-t}(L+1)\boldsymbol{v}^{(0)},T\right) \leq
CR^2\big\|\boldsymbol{v}^{(0)}\big\|_{H^4,m}^2=:
\overline{\mathcal{M}}.
\end{equation}
Now we will determine the threshold value by assuming that
$\boldsymbol{w}$ does not tend to zero. As a consequence of this assumption
we will derive the inequality
(\ref{eq:d14}) (equivalent to (\ref{eq:d15})
which bounds a norm of $\boldsymbol{v}^{(0)}$
from below. The inverse of this inequality gives the threshold which
we thus prove by contradiction.

After the above outline, we now turn to the details of the proof.
We thus assume that $\|\boldsymbol{w}(\cdot,T)\|_\infty$ does not
tend to zero as $T\rightarrow \infty$. According to
Theorem~\ref{thm:3} in Appendix~\ref{app:3} (with
$f=\|\boldsymbol{w}\|_{\tilde{H}_1}$) we then have
\[
\lim_{T\rightarrow \infty}\mathcal{N}(\boldsymbol{w},T)=\infty.
\]
In particular, there is a $T_0$ such that
\begin{equation}\label{eq:d11}
\mathcal{N}(\boldsymbol{w},T_0)=2\overline{\mathcal{M}}
\end{equation}
and
\begin{equation} \label{eq:d3}
\mathcal{N}(\boldsymbol{w},T)\leq 2\overline{\mathcal{M}},\qquad
T\leq T_0.
\end{equation}
Using (\ref{eq:d3}) and the Sobolev inequality (\ref{eq:a2}) we
can now estimate the three terms of
$\mathcal{M}(\boldsymbol{G}(\boldsymbol{w}),T_0)$ in terms of
$\overline{\mathcal{M}}$. For the first term we have the
inequality (\ref{eq:a11}) from Appendix~\ref{app:2},
\begin{equation} \label{eq:d6}
CR^2 \|\boldsymbol{G}(\boldsymbol{w})\|^2\leq CR^2
\|\boldsymbol{w}\|_{\infty}^2 \sum_{k=1}^3 \left\|\frac{\partial
\boldsymbol{w}}{\partial x_k}\right\|^2.
\end{equation}
The two factors in the right hand side can be estimate as follows
\begin{gather}
\|\boldsymbol{w}\|_{\infty}^2 \leq
C\|\boldsymbol{w}\|_{\tilde{H}_1}^2\leq C\overline{\mathcal{M}},
\label{eq:d4} \\
\left\|\frac{\partial \boldsymbol{w}}{\partial x_k}\right\|^2 \leq
C\|\boldsymbol{w}\|_{\tilde{H}_1}^2\leq C\overline{\mathcal{M}}.
\label{eq:d5}
\end{gather}
Inserting (\ref{eq:d4}) and (\ref{eq:d5}) into  (\ref{eq:d6}),
we obtain
\begin{equation} \label{eq:d7}
CR^2 \|\boldsymbol{G}(\boldsymbol{w})\|^2\leq
CR^2\overline{\mathcal{M}}^2.
\end{equation}
We use similar techniques and inequality (\ref{eq:a12}) and (\ref{eq:a13})
for the remaining two terms,
\begin{equation} \label{eq:d8}
CR^2\int_0^{T_0} \|\boldsymbol{G}(\boldsymbol{w})\|_{m}^2dt \leq
CR^2\|\boldsymbol{w}\|_{\infty}^2 \int_0^{T_0} \sum_{k=1}^3
\left\|\frac{\partial \boldsymbol{w}}{\partial
x_k}\right\|_{m}^2dt\leq CR^4\overline{\mathcal{M}}^2
\end{equation}
and
\begin{gather}
CR^2\int_0^{T_0} \|\boldsymbol{G}_t(\boldsymbol{w})\|_{m}^2dt \leq
CR^2 \sum_{k=1}^3 \left\|\frac{\partial \boldsymbol{w}}{\partial
x_k}\right\|_{m}^2\int_0^{T_0}\left\|\frac{\partial
\boldsymbol{w}}{\partial t}\right\|_{\infty}^2 dt \nonumber \\
\qquad {}+ CR^2\|\boldsymbol{w}\|_{\infty}^2\int_0^{T_0}
\sum_{k=1}^3\left\|\frac{\partial \boldsymbol{w}_t}{\partial
x_k}\right\|_{m}^2dt \leq CR^4\overline{\mathcal{M}}^2.
\label{eq:d10}
\end{gather}
Combining (\ref{eq:d7}), (\ref{eq:d8}) and (\ref{eq:d10})
we obtain
\begin{equation} \label{eq:d13}
\mathcal{M}(\boldsymbol{G}(\boldsymbol{w}),T_0) \leq
CR^4\overline{\mathcal{M}}^2.
\end{equation}
Now we insert (\ref{eq:d12}), (\ref{eq:d11}) and (\ref{eq:d13})
into (\ref{eq:d2})
to obtain
\begin{equation} \label{eq:d14}
2\overline{\mathcal{M}} \leq
\overline{\mathcal{M}}+CR^4\overline{\mathcal{M}}^2.
\end{equation}
We divide the inequality by $\overline{\mathcal{M}}$, insert the
expression (\ref{eq:d12}) for $\overline{\mathcal{M}}$ and
rearrange to
\begin{equation} \label{eq:d15}
\big\|\boldsymbol{v}^{(0)}\big\|_{H^4,m} ^2\geq \frac{1}{CR^6}.
\end{equation}
We have thus shown that an initial perturbation which leads to
instability must satis\-fy~(\ref{eq:d15}). If the initial
perturbation satisfies the inverse inequality then the instability
assumption leads to a contradiction. In other words, if the
initial perturbation satisfies the inverse inequality
of~(\ref{eq:d15}) then
$\|\boldsymbol{w}(\cdot,T)\|_{\infty}\rightarrow 0$ as $T$ tends
to infinity. Noting the $\boldsymbol{w}$ is related
to~$\boldsymbol{v}$ by~(\ref{eq:m10}), the proof of
Theorem~\ref{thm:1} is complete.

\appendix

\section{Definition and notation for norms} \label{app:1}

In this paper we track the $R$-dependence of all estimates we
derive. This forces us to introduce some non-standard norms which
weigh the different components of the vector function with
coefficients depending on $R$. Also we need a Sobolev inequality
(\ref{eq:a2}) with a~right hand side not containing all second
derivatives.

We start, however, to introduce notation for the standard norms.
The $L^2$-norm for functions $u:\Omega\rightarrow\mathbb{C}$ is defined by
\[
\|u\|^2=\int_\Omega |u(\boldsymbol{x})|^2 d\boldsymbol{x}.
\]
We use the same notation for the $L^2$-norm of vector functions
\[
\|\boldsymbol{u}\|^2=\sum_{k=1}^3\|u_k\|^2.
\]
It is clear from the argument if it is the scalar or vector norm which
is intended.
For the standard Sobolev spaces with ``integration exponent'' two we
use the following notation
\[
\|\boldsymbol{u}\|_{H^k}^2=\sum_{|\alpha|\leq k}
\left\|\frac{\partial^{|\alpha|}\boldsymbol{u}}{\partial
x^{\alpha}}\right\|^2,
\]
where $\alpha=(\alpha_1,\alpha_2,\alpha_3)$ is a multi-index.

We refer to the weighted norms also as modified norms and use the
subscript $m$ to identify them. Let $X$ denote any of the norms
introduced so far, then the corresponding modified version is
\[
\|\boldsymbol{u}\|_{X,m}^2=\|u_1\|_{X}^2+R^2\|u_2\|_{X}^2+\|u_3\|_{X}^2.
\]
In particular we have
\begin{equation} \label{eq:a4}
\|\boldsymbol{u}\|_{m}^2= \|u_1\|^2+R^2\|u_2\|^2+\|u_3\|^2.
\end{equation}
Before we give the Sobolev inequality (\ref{eq:a2}) we define the
semi-norm $J_2$ consisting of selected second order derivatives.
The choice of second order derivatives is dictated by the
difficulty in estimating normal derivatives (derivatives with
respect to the $x_2$-coordinate)
\[
J_2(\boldsymbol{u})=\left\|\frac{\partial^2\boldsymbol{u}}{\partial
x_1^2}\right\|^2+ \left\|\frac{\partial^2\boldsymbol{u}}{\partial
x_1
\partial x_2}\right\|^2+\left\|\frac{\partial^2
\boldsymbol{u}}{\partial x_2 \partial x_3}\right\|^2+
\left\|\frac{\partial^2 \boldsymbol{u}}{\partial x_3^2}\right\|^2.
\]
Now we can define two similar norms
\[
\|\boldsymbol{u}\|_{\tilde{H}_1}^2 =
\|\boldsymbol{u}\|_{m}^2+\sum_{k=1}^3 \left\|\frac{\partial
\boldsymbol{u}}{\partial x_k}\right\|^2+J_2(\boldsymbol{u})
\]
and
\[
\|\boldsymbol{u}\|_{\tilde{H}_2}^2 =
\|\boldsymbol{u}\|^2+\sum_{k=1}^3 \left\|\frac{\partial
\boldsymbol{u}}{\partial x_k}\right\|^2+J_2(\boldsymbol{u}).
\]
The $\tilde{H}_1$-norm is greater than the $\tilde{H}_2$-norm
($R\geq 1$) for which we have the Sobolev inequali\-ty
\begin{equation} \label{eq:a2}
\|\boldsymbol{u}\|_{\infty}^2 \leq C
\|\boldsymbol{u}\|_{\tilde{H}_2}^2,
\end{equation}
see \cite[p.~385]{KL:1}. By the left hand side in (\ref{eq:a2}) we
of course mean
\[
\|\boldsymbol{u}\|_{\infty}=\max_{k\in\{1,2,3\}} \mbox{ess}\,
\sup_{\boldsymbol{x}\in\Omega} |u_k(\boldsymbol{x})|.
\]

\section{Auxiliary estimates} \label{app:2}

To bound $\mathcal{M}_0(e^{-t}(L+1)\boldsymbol{v}^{(0)},T)$, we must
estimate the following four terms
\begin{gather*}
I=CR^2\big\|(L+1)\boldsymbol{u}^{(0)}\big\|^2,
\\
CR^2\int_0^T e^{-2t}\big\|(L+1)\boldsymbol{v}^{(0)}\big\|_m^2dt =
CR^2 \big\|(L+1)\boldsymbol{v}^{(0)}\big\|_m^2=:II,
\\
III=\big\|(L+1)\boldsymbol{v}^{(0)}\big\|_{\tilde{H}_2}^2
\end{gather*}
and
\[
IV=\big\|(L+1)^2\boldsymbol{v}^{(0)}\big\|_m^2.
\]
For the first term we have
\[
I\leq CR\big\|\Delta
\boldsymbol{v}^{(0)}\big\|^2+CR^2\left\|\frac{\partial
\boldsymbol{v}^{(0)}}{\partial x_1}\right\|^2
+CR^2\big\|v^{(0)}_2\big\|^2 \leq
CR^2\big\|\boldsymbol{v}^{(0)}\big\|_{H^2}^2.
\]
For the second term
\[
II\leq CR\big\|\Delta
\boldsymbol{v}^{(0)}\big\|_m^2+CR^2\left\|\frac{\partial
\boldsymbol{v}^{(0)}}{\partial x_1}\right\|_{m}^2
+CR^2\big\|v^{(0)}_2\big\|^2 \leq
CR^2\big\|\boldsymbol{v}^{(0)}\big\|_{H^2,m}^2.
\]
The third and fourth terms are bounded according to
\[
III\leq CR^2\big\|\boldsymbol{v}^{(0)}\big\|_{H^4,m}^2,\qquad
IV\leq CR^2\big\|\boldsymbol{v}^{(0)}\big\|_{H^4,m}^2.
\]
Combining the above results we obtain
\begin{equation} \label{eq:a10}
\mathcal{M}_0\left(e^{-t}(L+1)\boldsymbol{v}^{(0)},T\right)\leq
CR^2\big\|\boldsymbol{v}^{(0)}\big\|_{H^4,m}^2.
\end{equation}
To bound $\mathcal{M}(\boldsymbol{G}(\boldsymbol{w}),T)$ we need
the following estimates of the nonlinearity
\begin{gather}
\|\boldsymbol{G}(\boldsymbol{w})\|^2=\left\|\sum_{k=1}^3
w_k\frac{\partial \boldsymbol{w}}{\partial x_k}\right\|^2 \leq
\sum_{k=1}^3
\left\|w_k\frac{\partial \boldsymbol{w}}{\partial x_k}\right\|^2 \nonumber\\
\phantom{\|\boldsymbol{G}(\boldsymbol{w})\|^2}{}\leq \sum_{k=1}^3
\|w_k\|_{\infty}^2\left\|\frac{\partial \boldsymbol{w}}{\partial
x_k}\right\|^2 \leq \|\boldsymbol{w}\|_{\infty}^2\sum_{k=1}^3
\left\|\frac{\partial \boldsymbol{w}}{\partial
x_k}\right\|^2,\label{eq:a11} \\
\|\boldsymbol{G}(\boldsymbol{w})\|_{m}^2=\sum_{k=1}^3 \left(
\left\|w_k\frac{\partial w_1}{\partial x_k}\right\|^2+ R^2
\left\|w_k\frac{\partial w_2}{\partial x_k}\right\|^2+
\left\|w_k\frac{\partial w_3}{\partial x_k}\right\|^2
\right) \nonumber\\
\phantom{\|\boldsymbol{G}(\boldsymbol{w})\|_{m}^2}{}\leq
\sum_{k=1}^3 \|w_k\|_{\infty}^2\left\|\frac{\partial
\boldsymbol{w}}{\partial x_k}\right\|_{m}^2 \leq
\|\boldsymbol{w}\|_{\infty}^2 \sum_{k=1}^3\left\|\frac{\partial
\boldsymbol{w}}{\partial x_k}\right\|_{m}^2\label{eq:a12}
\end{gather}
and
\begin{gather}
\|\boldsymbol{G}_t(\boldsymbol{w})\|_{m}^2=\left\|\frac{\partial
}{\partial t}\left(\sum_{k=1}^3 w_k \frac{\partial
\boldsymbol{w}}{\partial x_k}\right)\right\|_{m}^2 \leq
\sum_{k=1}^3\left( \left\|\frac{\partial w_k}{\partial
t}\frac{\partial \boldsymbol{w}}{\partial x_k}\right\|_{m}^2+
\left\|w_k\frac{\partial \boldsymbol{w}_t}{\partial x_k}\right\|_{m}^2 \right)
\nonumber\\
 \phantom{\|\boldsymbol{G}_t(\boldsymbol{w})\|_{m}^2}{}\leq
\sum_{k=1}^3 \left\|\frac{\partial \boldsymbol{w}_k}{\partial
t}\right\|_{\infty}^2\left\|\frac{\partial
\boldsymbol{w}}{\partial x_k}\right\|_{m}^2+
\sum_{k=1}^3 \|w_k\|_{\infty}^2\left\|\frac{\partial \boldsymbol{w}_t}{\partial x_k}
\right\|_{m}^2 \nonumber\\
 \phantom{\|\boldsymbol{G}_t(\boldsymbol{w})\|_{m}^2}{}
\leq \left\|\frac{\partial \boldsymbol{w}}{\partial
t}\right\|_{\infty}^2\sum_{k=1}^3 \left\|\frac{\partial
\boldsymbol{w}}{\partial x_k}\right\|_{m}^2+
\|\boldsymbol{w}\|_{\infty}^2\sum_{k=1}^3 \left\|\frac{\partial
\boldsymbol{w}_t}{\partial x_k}\right\|_{m}^2. \label{eq:a13}
\end{gather}

\section{A Sobolev inequality for asymptotic stability} \label{app:3}

Let $f:[0,\infty)\rightarrow \mathbb{R}$. The theorem below gives
a condition which imply that $f(t)\rightarrow 0$ as the time ($t$)
tends to infinity. In this paper we will choose
$f(t)=\|\boldsymbol{w}(\cdot,t)\|_{\tilde{H}_1}$, combined with the inequality
(\ref{eq:a2}), the theorem can then be used to prove that
$\|\boldsymbol{w}(\cdot,t)\|_\infty\rightarrow 0$.

\begin{theorem} \label{thm:3}
If
\begin{equation} \label{eq:a18}
\int_0^\infty \left( |f(t)|^2+|f'(t)|^2 \right) dt<\infty,
\end{equation}
then
\[
\lim_{t\rightarrow \infty} f(t)=0.
\]
\end{theorem}

\begin{proof} The required Sobolev inequality is
\begin{equation} \label{eq:a19}
\sup_{t\geq T}|f(t)|^2 \leq
C\int_T^\infty \left( |f(t)|^2+|f'(t)|^2 \right) dt,
\end{equation}
where $C$ is independent of $T$. For a proof, see
\cite[appendix~3]{KL:1}.

Clearly
\begin{equation} \label{eq:a16}
|f(T)|^2\leq \sup_{t\geq T}|f(t)|^2.
\end{equation}
Because of (\ref{eq:a18}) we also have
\begin{equation} \label{eq:a17}
\lim_{T\rightarrow \infty}
\int_T^\infty \left( |f(t)|^2+|f'(t)|^2 \right) dt=0.
\end{equation}
Combining (\ref{eq:a19}), (\ref{eq:a16}) and (\ref{eq:a17}) completes the
proof.
\end{proof}

\label{Liefvendahl-lastpage}

\end{document}